\documentclass[11pt]{amsart}

\usepackage{amsmath}
\usepackage{amssymb}
\usepackage{mathtools}
\usepackage{amsthm}
\usepackage[foot]{amsaddr}
\usepackage{hyperref}
\usepackage{fullpage}
\usepackage{cleveref}
\usepackage{enumitem}
\usepackage{algpseudocode}
\usepackage{algorithm}
\usepackage{booktabs}
\usepackage{makecell}

\newcommand{\term}[1]{\textbf{#1}}
\newcommand{\tech}[1]{\textsf{#1}}

\newcommand{\zmat}[2]{\ensuremath{M_{#1 \times #2}(\Z)}}

\newcommand{\pdv}[2]{\ensuremath{\frac{\partial #1}{\partial #2}}}
\newcommand{\dv}[2]{\ensuremath{\frac{\mathrm{d}\, #1}{\mathrm{d}\, #2}}}

\newcommand{\bolda}{\ensuremath{\mathbf{a}}}

\newcommand{\boldc}{\ensuremath{\mathbf{c}}}

\newcommand{\bolde}{\ensuremath{\mathbf{e}}}

\newcommand{\boldx}{\ensuremath{\mathbf{x}}}
\newcommand{\boldy}{\ensuremath{\mathbf{y}}}
\newcommand{\boldz}{\ensuremath{\mathbf{z}}}

\newcommand{\boldone}{\boldsymbol{1}}
\newcommand{\boldzero}{\boldsymbol{0}}

\newcommand{\Z}{\mathbb{Z}}

\newcommand{\Q}{\mathbb{Q}}
\newcommand{\C}{\mathbb{C}}

\newcommand{\euler}{  \mathcal{E} }
\newcommand{\newton}{ \mathcal{N} }
\newcommand{\inner}[2]{ \ensuremath{\left\langle {#1}, {#2}\right\rangle} }
\newcommand{\homog}[1]{ #1^{\operatorname{hom}} }

\DeclareMathOperator{\diag}{diag}

\theoremstyle{plain}
\newtheorem{theorem}{Theorem}

\theoremstyle{definition}

\theoremstyle{definition}
\newtheorem{remark}[theorem]{Remark}

\title{GPU-accelerated path tracker for polyhedral homotopy}
\author{Tianran Chen}
\address{Department of Mathematics, Auburn University at Montgomery, Montgomery Alabama USA}
\email{ti@nranchen.org}

\begin{document}

\begin{abstract}
    The polyhedral homotopy method of Huber and Sturmfels
    is a particularly efficient and robust numerical method
    for solving system of (Laurent) polynomial equations.
    A central component in an implementation of this method
    is an efficient and scalable path tracker.
    While the implementation issues in a scalable path tracker for 
    computer clusters or multi-core CPUs have been solved thoroughly,
    designing good GPU-based implementations is still an active research topic.
    This paper addresses the core issue of 
    efficiently evaluate a multivariate system of Laurent polynomials
    together with all its partial derivatives.
    We propose a simple approach that maps particularly well onto the
    parallel computing architectures of modern GPUs.
    As a by-product, we also simplify and accelerate the
    path tracker by consolidating the computation of Euler and Newton directions.
\end{abstract}

\maketitle


\section{Introduction}

Numerical homotopy continuation methods have emerged as a robust and efficient
family of numerical methods for solving systems of polynomial equations
\cite{BatesHauensteinSommeseWampler2013Numerically,ChenLi2015Homotopy,sommese_numerical_2005}.
Their most salient advantage is that each solution can be computed independently
making such methods \emph{pleasantly parallel}.
A great variety of specific homotopy constructions have been proposed.
In their seminal work \cite{HuberSturmfels1995Polyhedral},
Huber and Sturmfels introduced the particularly attractive
\emph{polyhedral homotopy} method,
which can take full advantage of
the Newton polytope structure of the target system.
The implementation of polyhedral homotopy on computer clusters and multi-core CPUs
has been explored thoroughly.
This paper addresses the unique challenges in leveraging the power of GPUs
and similar devices.

Central to the polyhedral homotopy method (and homotopy methods in general)
are the ``path tracking'' algorithms for tracking the homotopy paths.
While a great variety of numerical methods can be used for path tracking,
the \emph{predictor-corrector} scheme has emerged as the method of choice
within the community of numerical homotopy methods due to its superior efficiency
and stability \cite{BatesHauensteinSommeseWampler2013Numerically,Li2003Numerical,SommeseWampler2005Numerical}.
For a typical predictor-corrector based path tracker,
a significant portion of computation time is devoted to the
evaluation of the homotopy function together with its partial derivatives.
Indeed, this part may dominate the overall computational cost
in the path tracking process
for systems having complicated expressions relative to its dimension.
The simultaneous evaluation of a (multivariate) system of polynomial functions
given by the polyhedral homotopy construction together with 
all its partial derivatives is our main focus.
Tremendous efforts have been devoted to this problem within the community of
polyhedral homotopy method
(e.g. \cite{Kojima2008,LeeLiTsai2008HOM4PS-2.0,VerscheldeYoffe2012evaluating,VerscheldeYu2016Polynomial}) 
as well as the community of the homotopy method in general.
In the broader context, this is the main subject of the established field of
\emph{Automatic Differentiation} \cite{Griewank2008Evaluating}.

In the present contribution we propose a very simple approach for
efficiently evaluating a homotopy function together with all its partial derivatives,
in the context the polyhedral homotopy method,
that is particularly efficient in modern GPUs and similar devices.
Indeed, the problem of evaluating Laurent polynomial systems
together with their derivatives is converted into 
general matrix multiplication (\tech{GEMM}) operations \cite{BLAS},
which modern GPUs
(especially the newer designs like the ``Tensor Cores'' of \tech{NVIDIA})
are designed to optimize.
Grounded on this technique, we also unify the computation of
Euler and Newton directions that are needed in path tracking.
\section{A review of polyhedral homotopy} \label{sec: polyhedral}

This section briefly review the polyhedral homotopy of Huber and Sturmfels
\cite{HuberSturmfels1995Polyhedral}. 
For simplicity, we first restrict our attention to systems having 
\emph{``generic coefficients''} in the sense of the Bernshtein's First Theorem 
\cite{Bernshtein1975Number}.
A square \emph{Laurent polynomial system} in the variables
$\boldx = (x_1,\dots,x_n)$ is a system $F = (f_1,\ldots,f_n)$ given by
$
    f_k(\boldx) = \sum_{\bolda\in S_k} c_{1,\bolda} \, \boldx^{\bolda}
$
for $k=1,\ldots,n$,
where the nonempty finite set of column vectors $S_k \subset \Z^n$
describes the exponents appeared in the monomials in $f_k$,
and the ``multi-exponent'' notation
$
    \boldx^{\bolda} = 
    x_1^{a_1} \, \cdots \, x_n^{a_n}
$
describes a monomial whose exponents are given by the (column) vector $\bolda$.
The goal is to find \emph{all} isolated zeros of $F(\boldx)$
in $(\C^*)^n = (\C \setminus \{ \boldzero \})^n$.
Given \emph{lifting functions} $\omega_k : S_k \to \Q^+$ for $k=1,\dots,n$, 
the \emph{polyhedral homotopy} of Huber and Sturmfels
for solving the target system $F(\boldx) = \boldzero$
is given by the homotopy function $H = (h_1,\ldots,h_n)$ with
\begin{equation}
    \label{equ:polyhedral}
    h_k(\boldx,t) =
    \sum_{\bolda\in S_k} 
        c_{k,\bolda} \,
        \boldx^{\bolda} \,
        t^{\omega_k(\bolda)}
    \quad
    \text{for } k = 1, \ldots, n.
\end{equation}
Clearly, 
$H(\boldx,1) \equiv F(\boldx)$ is exactly the target system.
Under the genericity assumption,
$H(\boldx,t)=\boldzero$ defines smooth paths
in $(\C^*)^n \times (0,1]$, parametrized by $t$,
emanating from isolated $\C^*$-zeros of 
$F(\boldx) \equiv H(\boldx,1)$ and continue toward $t=0$.
Once the starting points of these paths at $t=0$ are obtained,
the $\C^*$-zeros of the target system $F(\boldx)$
can be found by tracking these paths.

Since the starting points are not directly involved
in our discussions, we simply assume 
the set of all solution to $H(\boldx,t_0) = \boldzero$
in $(\C^*)^n$ for a common $t$-value of $t=t_0 > 0$,
is readily available and refer to
Refs.~\cite{HuberSturmfels1995Polyhedral}
for detailed discussions on the step of locating the starting points.

\subsection{Logarithmic homotopy parameter}

One potential source of numerical instability
is the scale of $\| \pdv{H}{t} \|$ near $t=0$, which may be infinite. 
The most widely accepted solution
is to use the substitution
$t = e^\tau$
and adopt $\tau$ as the path parameter \cite{gunji_phompolyhedral_2004,LeeLiTsai2008HOM4PS-2.0}.
The original interval $t \in (0,1)$ is maps to 
the corresponding interval $\tau \in (-\infty,0)$,
thus 
``dilutes'' the numerical instability.
With this, \eqref{equ:polyhedral} can be reformulated as
$H = (h_1,\ldots,h_n) : (\C^*)^n \times (-\infty,0]$ given by
\begin{equation}\label{equ:polyhedral-exp-t}
    h_k(\boldx,\tau) =
    \sum_{\bolda\in S_k} 
        c_{k,\bolda} \,
        \boldx^{\bolda} \,
        e^{\tau \cdot \omega_k(\bolda)}.
\end{equation}
In practice, the path tracking process start from some $t_0 < 0$
with sufficiently large magnitude.
We refer to Refs.~\cite{gunji_phompolyhedral_2004,LeeLiTsai2008HOM4PS-2.0}
for detailed descriptions. 
The downside%
\footnote{
    Another apparent downside is the use of
    potentially expensive exponential function.
    However, the original formulation \eqref{equ:polyhedral}
    already requires the computation of rational powers
    $t^{\omega_k(\bolda)}$, which are computed indirectly 
    via exponential functions,
    so the use of exponential functions is unavoidable
    in most practical implementations.
}
is the potentially longer path parameter interval $[t_0,0]$.
However, this technique is adopted by several polyhedral homotopy implementations
\cite{ChenLeeLi2014Hom4PS-3,gunji_phompolyhedral_2004,LeeLiTsai2008HOM4PS-2.0}
and the consensus is that the numerical benefits far outweighs the cost.

\subsection{Projective formulation}

Another issue one has to deal with is the existence of divergent paths,
especially paths along which $\|\boldx\| \to \infty$ or grow very large.
This is usually dealt with by lifting \eqref{equ:polyhedral-exp-t}
into the complex projective space $\mathbb{P}^n$.
Using $\boldy = ( y_0 , y_1 , \cdots , y_n )$
that represents homogeneous coordinates in $\mathbb{P}^n$,
\eqref{equ:polyhedral-exp-t} may be further reformulated as
$\homog{H} : \C^{n+1} \times (-\infty,0]$ with
\begin{equation}\label{equ:polyhedral-hom}
    \homog{h}_k(\boldy,\tau) =
    \sum_{\bolda \in S_k} 
        c_{k,\bolda} \,
        \boldy^{\homog{\bolda}} \,
        e^{\tau \cdot \omega_k(\bolda)},
        \quad\text{for } k=1,\ldots,n,
\end{equation}
where 
$\homog{\bolda} = (\bolda , \deg(f_i) - \boldone^\top \bolda)^\top \in \Z^{n+1}$,
for each $\bolda \in S_k \subset \Z^n$,
represent the exponent vectors of the corresponding monomials
in the homogenization of $f_k$.
Then $\homog{H} = (\homog{h}_1, \ldots, \homog{h}_n)$
is homogeneous in $\boldy$,
and represents the lifting of the homotopy \eqref{equ:polyhedral-exp-t}
into the complex projective space.
We refer to Refs.~\cite{BlumRivest1988training,ChenLi2012Spherical,ShubSmale1993Complexity,SommeseWampler2005Numerical}
for discussions on this topic.

\subsection{Euler and Newton directions in path tracking}

In the affine formulation \eqref{equ:polyhedral-exp-t},
the paths in $(\C^*)^n \times (-\infty,0]$ are defined by
$H(\boldx,\tau) = \boldzero$.
Along each, by the Implicit Function Theorem, 
$\boldx$ is a smooth function of $\tau$,
and $\frac{d \boldx}{d\tau}$, as a column vector, is characterized by
the Davidenko equation
\begin{align*}
    \frac{\partial H}{\partial \boldx}(\boldx,\tau)
    \frac{d \boldx}{d \tau}(\boldx,\tau) +
    \frac{\partial H}{\partial \tau}(\boldx,\tau)
    &= \boldzero.
\end{align*}
This vector shall be referred to as the \emph{Euler direction}
in affine coordinates of a solution path at a given point on the path,
and we shall use the notation
$\euler(\boldx,\tau) := \frac{d \boldx}{d\tau}(\boldx,\tau)$.

The corresponding \emph{Euler direction}
$\euler(\boldy,\tau)$ in homogeneous coordinates is characterized by
\begin{align}\label{equ: euler direction}
    \frac{\partial \homog{H}}{\partial \boldy}(\boldy,\tau)
    \euler(\boldy,\tau) +
    \frac{\partial \homog{H}}{\partial \tau}(\boldy,\tau)
    &= \boldzero,
    &&\text{and} &
    \boldy^* \, \euler(\boldy,\tau) &= 0,
\end{align}
where $\boldy^*$ is the row vector whose entries are
the conjugates of the entries of $\boldy$.
$\euler(\boldy,\tau)$ represents the geodesic direction 
defined by the path with respect to the
Fubini-Study metric on $\mathbb{P}^n$.
We again refer to Ref.~\cite{ChenLi2012Spherical} for the derivation.

The name, ``Euler direction'' refers to the fact that in both form,
$\Delta\tau \mapsto \boldx(\tau) + \Delta \tau \euler(\boldx,\tau)$ and
$\Delta\tau \mapsto \boldy(\tau) + \Delta \tau \euler(\boldy,\tau)$ are
the first order approximations of the solution path
provided by the Euler's method,
in the affine and homogeneous coordinates, respectively.
Under the smoothness assumption,
both are uniquely defined.
They can be used in any predictor methods
that utilize first order derivative information.
In particular, they are directly used in
Euler's method, Cubic Hermite method, and certain spline-based methods.
They are also indirectly used in predictors based on
Runge-Kutta (Runge-Kutta-Fehlberg) methods.

Similarly, the \emph{Newton direction},
in affine coordinates, is the unique vector $\newton(\boldx,\tau)$
such that
\begin{equation*}
    \frac{\partial H}{\partial \boldx}(\boldx,\tau) \newton(\boldx,\tau) + 
    H(\boldx,\tau)
    = \boldzero,
\end{equation*}
and the corresponding \emph{Newton direction} in homogeneous coordinates,
$\newton(\boldy,\tau)$ is defined by
\begin{align}\label{equ: newton direction}
    \frac{\partial \homog{H}}{\partial \boldy}(\boldy,\tau) \newton(\boldy,\tau) + 
    \homog{H}(\boldy,\tau)
    &= \boldzero.
    &&\text{and} &
    \boldy^* \, \newton(\boldy,\tau) &= 0.
\end{align}
They are vectors used in standard and projective Newton
\cite{BlumRivest1988training,ShubSmale1993Complexity}
iterations respectively.
The affine version can also be used in 
dampened Newton iterations
and local Newton homotopy method.%

\section{Hardware architectures and the problem statement}\label{sec: problem}

Homotopy continuation methods, including the polyhedral homotopy described above, 
are naturally parallel and scalable
in the sense that every solution can be computed independently
(pleasantly parallel).
The problem of implementing such parallel homotopy-based solvers on
traditional computer clusters and multi-core CPU architectures
has been thoroughly explored.
However, today, the landscape of high performance computing
is dominated by GPUs (graphical processing units),
and the same type of algorithms cannot be ported directly to modern GPUs
with the same level of efficiency and scalability
due to the inherent hardware constraints:
At the lowest level, modern GPUs gain their superiority in computational power
from their single-instruction-multiple-data (or \tech{SIMD}) architecture,
in which multiple processing unit perform the exact same mathematical operation
on multiple data points simultaneously.
Newest CPUs also support \tech{SIMD} models
(e.g., through \tech{Intel}'s \tech{AVX} instructions) to a certain extend.
Such architectures exploit data-parallelism, but not true concurrency,
hence the limitations.


Single-instruction-multiple-threads
(or \tech{SIMT}) is a subclass of \tech{SIMD} in
Flynn's 1972 taxonomy \cite{Flynn1972Computer}.
It is loosely defined by the parallel computing units having
their own independent registers and memory (cache and data Memory),
which better describes the higher level data organization in modern GPUs.
We will also include the structurally similar Tensor Cores%
\footnote{
    \tech{NVIDIA} introduced dedicated Tensor Cores
    for general matrix multiplications (\tech{GEMM})
    in 2017 with \tech{Volta V100} GPUs.
    \tech{Google}'s \tech{TPU}s, or Tensor Processing Units,
    have similar design and purpose.
}
in our target model.
Tensor cores specialize in general matrix multiplications (\tech{GEMM}),
which also resembles a \tech{SIMT} model from an algorithmic point of view.
To target these hardware devices, we assume a heterogeneous memory model
in which threads have access to different memory types 
with different latency and bandwidth expectations.%
\footnote{
    In \tech{NVIDIA} GPUs, for example,
    the different types of memory include
    host memory (including pinned host memory),
    global memory, and shared memory,
    as well as constant and texture memory that
    not directly relevant here.
    In addition, each thread also has access to its own
    register file and local (host) memory.
    The differences in the latency and bandwidth in
    accessing these memory types span several orders of magnitude,
    and the usage of these memory types must be considered carefully.
}





For simplicity, the hardware architectures with all these constraints,
including  
GPUs, CPU with \tech{AVX} instruction sets, and the new tensor cores
will simply be referred to as the \tech{SIMD/SIMT} architecture.
The computational framework proposed in this paper
directly targets such a \tech{SIMD/SIMT} architecture
and only this architecture.

The central problem here is to\ldots
\begin{enumerate}
    \item 
        evaluate $H$ and its Jacobian matrix 
        (i.e., all its partial derivatives) simultaneously, and
    \item 
        compute the Euler and Newton directions simultaneously,
\end{enumerate}
while optimizing the algorithm subject to the constraints of
the \tech{SIMD/SIMT} hardware architecture.
In particular, the algorithm must\ldots
\begin{enumerate}[label=(\alph*)]
    \item minimize data transfer,
    \item minimize global memory access,
    \item minimize communication between devices, and
    \item extract parallelism from multiple levels
        (e.g., monomial, variable, path, and hardware level).
\end{enumerate}
In \Cref{sec:eval,sec: euler newton},
we propose a simple algorithm that solve this problem
under the given hardware constraints.

\section{Related works} \label{sec: related}

Evaluating multivariate systems together with their Jacobian matrices
efficiently is a basic problem in numerical computation
(see, e.g., Ref.~\cite{CarnicerGasca1990Evaluation}).
In the context of homotopy continuation methods,
this problem was analyzed rigorously in Ref.~\cite{Kojima2008}
through the framework of multivariate generalization of the Horner's rule.
Several algorithms have been designed specifically for GPUs
\cite{VerscheldeYoffe2012evaluating,VerscheldeYu2016Polynomial}.
Several algorithms originally designed for CPUs
\cite{ChenLeeLi2014Hom4PS-3,LeeLiTsai2008HOM4PS-2.0}
have also been ported to \tech{CUDA}-based GPUs.

One distinguishing feature of the algorithm proposed in this paper
is its stark simplicity.
More importantly, the proposed approach formulates the central problem
into problems that GPUs are optimized --- ``\tech{BLAS} operations''
\cite{BLAS_updated,BLAS}.  

\section{Formulating homotopy evaluation as matrix multiplication} \label{sec:eval}

\subsection{Logarithmic formulation}

For simplicity, we assume the target system $F = (f_1,\ldots,f_n)$
to be ``unmixed'', i.e., the support of $f_1,\dots,f_n$ are identical.
We also restrict ourselves to solutions within $(\C^*)^n$
like the original polyhedral homotopy formulation.
Under these assumptions, we formulate the problem of evaluating
the homotopy function $\homog{H}$, in \eqref{equ:polyhedral-hom}
together with all its partial derivatives
as a problem of general matrix multiplication (\tech{GEMM}),
which maps nice to the \textsf{SIMD/SIMT} architectures.

We first lift $H$ into the logarithmic coordinate space.
Fixing any branch $\log : \C^* \to \C$ of the complex logarithm function, 
we introduce new variables $\boldz = (z_0,\dots,z_n)$
through the relations
\begin{equation}\label{equ:z-log-x}
    \boldz = \log \boldy = (\log y_0, \dots, \log y_n).
\end{equation}
Though this depends on the choice of the branch of $\log$, 
the projection $\boldy = e^{\boldz} = (e^{z_1},\dots,e^{z_n})$
remains valid and well-defined 
regardless of the choices.
Thus, under the restriction $\boldy \in (\C^*)^n$,
this represents a locally smooth change of variables.
Under the ``unmixed'' assumption,
we define $S := \{ \bolda_1, \ldots, \bolda_m \} = S_1 = \cdots = S_n$
to be the common support of $f_1,\ldots,f_n$,
and let $\{ \omega_1, \ldots, \omega_m \}$ be the corresponding lifting values.
We also use the notations
\begin{align*}
    N &= n + 1 \\
    \hat{\boldz} &= (\boldz,\tau), \\
    \homog{\hat{\bolda}}_k &= (\homog{\bolda}_k, \omega_k)^\top \quad\text{for } k=1,\dots,m, \text{ and }\\
    \homog{\hat{A}} &= 
    \begin{bmatrix}
        \homog{\hat{\bolda}}_1 & \cdots & \homog{\hat{\bolda}}_m
    \end{bmatrix}
    \in \zmat{(N+1)}{m}.
\end{align*}
Then the collection of monomials in $\homog{h}_k(\boldy,\tau)$,
as a row vector, can be expressed as
\begin{equation}
    \begin{bmatrix}
        \boldy^{\homog{\bolda}_1} e^{\tau \cdot \omega_1} 
        & \cdots &
        \boldy^{\homog{\bolda}_m} e^{\tau \cdot \omega_m}
    \end{bmatrix} =
    \begin{bmatrix}
        e^{\inner{\hat{\boldz}}{\homog{\hat{\bolda}}_1}}
        & \cdots &
        e^{\inner{\hat{\boldz}}{\homog{\hat{\bolda}}_m}}
    \end{bmatrix} =
    e^{\hat{\boldz} \homog{\hat{A}}}
\end{equation}
This leads to the compact expression
\begin{equation}
    \label{equ:exp-h}
    \homog{h}_k(\boldy,\tau) = 
    e^{\hat{\boldz} \homog{\hat{A}}} \cdot \boldc_k
    \quad\text{for } k = 1,\dots,n
\end{equation}
where $\boldc_k = (c_{k,1},\dots,c_{k,m})^\top$ is the column vector
that collects the coefficients of the terms in $\homog{h}_k$ for $k=1,\ldots,n$.

The logarithmic coordinate turns the computation of monomials into 
a simple vector-matrix product $\hat{\boldz} \mapsto \hat{\boldz} \homog{\hat{A}}$,
which is particularly efficient to compute on \tech{SIMD/SIMT} architectures.
However, the true benefit of the formulation \eqref{equ:exp-h} is the ease
with which partial derivatives can be computed using existing data. 
In particular,
\[
    \frac{\partial \homog{h}_k}{\partial z_j} =
    \frac{\partial}{\partial z_j}
    \left( e^{\hat{\boldz} \homog{\hat{A}}} \cdot \boldc_k \right)
    =
    e^{\hat{\boldz} \homog{\hat{A}}} \cdot 
    \left[
    \begin{smallmatrix}
        c_{k,1} \homog{a_{j,1}} \\
        \vdots \\[5pt]
        c_{k,m} \homog{a_{j,m}}
    \end{smallmatrix}
    \right]
    \quad\text{for } j = 0,\ldots,n.
\]
Notice that the vector appeared at the end is precisely the entrywise product
(Hadamard product) between the coefficient vector $\boldc_k$
and the transpose of the $j$-th row of $\homog{\hat{A}}$. 
To further simplify the notation as well as the data organization,
we define $(N+2) \times m$ matrices
\begin{equation*}
	B_k = 
	\begin{bmatrix}
		c_{k,1} \homog{\hat{\bolda}}_1 & \cdots & c_{k,m} \homog{\hat{\bolda}}_m \\
		c_{k,1}                        & \cdots & c_{k,m}
	\end{bmatrix},
    \quad\text{for } k=1,\dots,n,
\end{equation*}
which will remain constant and thus can be pre-computed
prior to the path tracking process.
Then, for any $k=1,\dots,n$, the function value $\homog{h}_k$
and all its partial derivatives with respect to 
$(z_0,\dots,z_n,\tau)$ and 
$(y_0,\dots,y_n,\tau)$ 
can be computed as the row vectors
\begin{align*}
    \begin{bmatrix}
	    D_{\boldz} \homog{h}_k & D_{\tau} \homog{h}_k & \homog{h}_k
    \end{bmatrix}
    &=
    \begin{bmatrix}
        \frac{\partial \homog{h}_k}{\partial z_0}
        & \cdots &
        \frac{\partial \homog{h}_k}{\partial z_n} &
        \frac{\partial \homog{h}_k}{\partial \tau} &
        \homog{h}_k
    \end{bmatrix} =
    e^{\hat{\boldz} \homog{\hat{A}}} B_k^\top, \\
    \begin{bmatrix}
	    D_{\boldx} \homog{h_k} & D_{\tau} \homog{h_k} & \homog{h_k}
    \end{bmatrix}
    &=
    \begin{bmatrix}
        \frac{\partial \homog{h_k}}{\partial y_0}
        & \cdots &
        \frac{\partial \homog{h_k}}{\partial y_n} &
        \frac{\partial \homog{h_k}}{\partial \tau} &
        \homog{h}_k
    \end{bmatrix} =
    e^{\hat{\boldz} \homog{\hat{A}}} B_k^\top 
    \diag( e^{-\boldz}, 1, 1 ).
\end{align*}

To apply these formulations to all components of $\homog{H}$,
it is convenient to define the \term{extended Jacobian matrices}
in logarithmic and homogeneous coordinates to be
\begin{align*}
    \hat{D}_{\hat{\boldz}} \homog{H} &:= 
    \begin{bmatrix} D_\boldz \homog{H} & D_\tau \homog{H} & \homog{H} \end{bmatrix}
    &&\text{and} &
    \hat{D}_{\hat{\boldy}} \homog{H} &:= 
    \begin{bmatrix} D_\boldy \homog{H} & D_\tau \homog{H} & \homog{H} \end{bmatrix}.
\end{align*}
Then the vectorization of these matrice can be computed as
\begin{align}
    \operatorname{vec}(\hat{D}_{\hat{\boldz}} \homog{H})
    &=
    e^{\hat{\boldz} \hat{A}} 
    \begin{bmatrix}
        B_1^\top & \cdots & B_n^\top
    \end{bmatrix}
    \\
    \operatorname{vec}(\hat{D}_{\hat{\boldy}} \homog{H})
    &=
    e^{\hat{\boldz} \hat{A}} 
    \begin{bmatrix}
        B_1^\top \diag( e^{-\boldz}, 1, 1 ) 
        & \cdots &
        B_n^\top \diag( e^{-\boldz}, 1, 1 )
    \end{bmatrix}.
\end{align}

\subsection{Batched evaluation at a large number of points}

The main advantage of \textsf{SIMD/SIMT} architectures,
such as GPU devices, is their ability in carrying out
identical and independent operations
on large sets of data at the same time.
Indeed, operating on sufficiently large data set
is the main mechanism by which latency can be hidden.
This characteristic maps well to the path tracking problem
due to its ``pleasantly parallel'' nature.

For $p$ different points $\hat{\boldz}_1,\dots,\hat{\boldz}_p$
in $\C^{n+1}$ with $\hat{\boldz}_k = (\boldz_k, \tau_k)$,
the vectorization of the extended Jacobian matrix
$\hat{D}_{\hat{\boldz}} \homog{H}$ evaluated at these points
can thus be computed simultaneously as
\begin{align*}
    \begin{bmatrix}
        \operatorname{vec}(\hat{D}_{\hat{\boldz}} \homog{H} (\hat{\boldz}_1)) \\
        \vdots \\
        \operatorname{vec}(\hat{D}_{\hat{\boldz}} \homog{H} (\hat{\boldz}_p)) \\
    \end{bmatrix}
    &=
    e^{\hat{Z} \hat{A}}
    \begin{bmatrix}
        B_1^\top & \cdots & B_n^\top
    \end{bmatrix}
    &
    &\text{where}
    &
    \hat{Z} &=
    \begin{bmatrix}
        \hat{\boldz}_1 \\
        \vdots \\
        \hat{\boldz}_p
    \end{bmatrix},
\end{align*}
and, with a slight abuse of notations,
we consider $e^{\hat{Z} \hat{A}}$ to be a $p \times m$ matrix
whose $k$-th row is $e^{\hat{z}_k \hat{A}}$.
Its counterpart in homogeneous coordinate can be computed
from the above formula by multiplying
$\diag(e^{-\boldz_i},1,1)$ to each block.

The simultaneous evaluation of $\homog{H}$ and all its partial derivatives,
i.e., the extended Jacobian matrix,
at a large number of points is thus formulated as a
general matrix multiplication (\tech{GEMM}) operation
(a level-3 \tech{BLAS} operaton),
which is one of the operations GPUs (especially tensor cores)
are optimized to carry out.
When the number of points $p$ is sufficiently large,
this operation can hide latency effectively.



\section{Simultaneous computation of the Euler and the Newton direction}
\label{sec: euler newton}

We extract an additional level of parallelism
from the simultaneous computation of the Euler and Newton direction
through a single unified process.
Even though the two directions are rarely used together,
this design consolidates the two independent steps,
which eliminates the costs associated with synchronizations
and greatly simplifies the data pathways.
More importantly, it eliminates complicated branching
(i.e., high level decision-making) from the GPU side.

To take full advantage of the \tech{SIMD/SIMT} architecture,
one must leverage of the benefit of ``batched'' operations,
i.e., structurally identical matrix operations
applied to a large number of matrices of the same size
and memory layout.
Not all matrix operations can be consolidated into batch mode,
and thus the operations must be chosen with care..

Recall that the Euler direction $\euler(\boldy,\tau)$ and 
the Newton direction $\newton(\boldy,\tau)$ 
are defined by the nonsingular linear systems
\begin{align*}
    \frac{\partial \homog{H}}{\partial \boldy} \euler(\boldy,\tau) +
    \frac{\partial \homog{H}}{\partial \tau}
    &= \boldzero,
    &&\text{and} &
    \boldy^* \, \frac{d \boldy}{d \tau} &= 0; \\
    \frac{\partial \homog{H}}{\partial \boldy} \newton(\boldy,\tau) + \homog{H}
    &= \boldzero.
    &&\text{and} &
    \boldy^* \, \Delta \boldy &= 0.
\end{align*}


The similarity between the two suggests that the equations
can be unified into the system
\begin{equation}
    \begin{bmatrix}
        \pdv{\homog{H}}{\boldy} & \pdv{\homog{H}}{\tau} & H \\
        \boldy^*        & 0          & 0
    \end{bmatrix}
    \begin{bmatrix}
        \dv{\boldy}{\tau} & \Delta \boldy \\
        1 & 0 \\
        0 & 1
    \end{bmatrix}
    =
    \mathbf{0}.
\end{equation}
Note that the matrix on the left contains the extended Jacobian matrix
$\hat{D}_{\hat{\boldy}} \homog{H}$.
Therefore, they can be computed through a numerically stable QR decomposition.
After finding an $(N+2) \times (N+2)$ unitary matrix $Q$ 
and an $N \times N$ upper triangular matrix $R$ such that
\[
    J^\top :=
    \begin{bmatrix}
        \pdv{\homog{H}}{\boldy} & \pdv{\homog{H}}{\tau} & H \\
        \boldy^*        & 0          & 0
    \end{bmatrix}^\top
    = Q \begin{bmatrix} R \\ \boldzero_{2 \times n} \end{bmatrix},
\]
then the null space of $J$
is spanned by the two vectors 
$(\bolde_{N+1} Q^*)^\top$ and $(\bolde_{N+2} Q^*)^\top$,
i.e., the complex conjugate of the right most two columns of $Q$.
After further reduction, the Euler and Newton directions are computed.

If we further assume $\hat{D}_{\hat{\boldy}} \homog{H}$ to be well conditioned,
then this process can be carried out via simple Householder transformations,
which, on most modern \tech{SIMD/SIMT} architectures,
can be process in ``batch mode'' on a large number of points simultaneously.

\begin{remark}
    In particular, on \tech{NVIDIA CUDA} architecture,
    this can be achieved through ``batched QR-decomposition''
    offered by \tech{cuBLAS} library.
\end{remark}

\begin{remark}
    If $\hat{D}_{\hat{\boldy}} \homog{H}$ is ill-conditioned,
    more advanced QR decomposition methods have to be used.
    Currently, these are not supported by the standard \tech{NVIDIA CUDA} framework.
    One potential solution is to use more advance solver libraries such as \tech{MAGMA}.
    Alternatively, pre-conditioners can be used to improve the numerical conditions
    of $\hat{D}_{\hat{\boldy}} \homog{H}$ prior to this step.
    The pros and cons of these solutions will be investigated in future studies.
\end{remark}

\section{Summary of the algorithm}

In the following, we summarize the main algorithms
and provide pseudocode for an implementation targeting GPUs.
We assume the target GPU has heterogeneous memory organization,
i.e., there is a distinction between ``host'' (CPU) memory
and ``device'' (GPU) memory, and transferring data between the two
is slow and must be minimized.
Operations in pseudocode will be marked by
``(Host)'', ``(Device)'', or ``(Host$\to$device)''
to indicate if an operation access host memory, device memory,
or transfer data from host memory to device memory.
There is also a nonhomogeneous hierarchy of different device memory types.
Their distinctions are left as implementation details
and not directly described here.

\Cref{alg: initialize} allocate device memory and initialize
constant matrices to be used in later computations.
Here $p \in \Z^n$ is the maximum number of points, i.e., paths,
to be stored simultaneously on a GPU device.
This number is mainly limited by the available global memory on device.

\begin{algorithm}
    \caption{Initialization procedure}
    \label{alg: initialize}
    \begin{algorithmic}[1]
        \Procedure{Initialize}{$p,\boldc_1,\ldots,\boldc_n,\homog{\hat{A}}$}
            \State (Host$\to$device) Transfer coefficient vectors $\boldc_1,\ldots,\boldc_n$ to device;
            \State (Host$\to$device) Transfer exponent matrix $\homog{\hat{A}}$ to device;
            \State (Device) Allocate $p \times (N+1)$ $\C$-matrix for $\hat{Z}$;
            \State (Device) Allocate $p \times (N+1)$ $\C$-matrix for $\hat{Y}$;
            \State (Device) Allocate $p(n+1) \times (N+2)$ $\C$-matrix for $J$;
            \State (Device) Allocate $(N+1) \times m$ $\C$-matrix for $B_k$ for each $k=1,\ldots,n$;
            \State (Device) Allocate $2p \times (N+2)$ $\C$-matrix for $V_i$ for each $i=1,\ldots,n$;
            \State (Device) $B_k \gets \begin{bmatrix} \homog{\hat{A}} \diag(\boldc_k) \\ \boldc_k \end{bmatrix}$ for $k=1,\ldots,n$;
        \EndProcedure
    \end{algorithmic}
\end{algorithm}

\begin{algorithm}
    \caption{Batched evaluation of extended Jacobian matrix}
    \label{alg: evaluation}
    \begin{algorithmic}[1]
        \Procedure{Evaluate}{$i,b$}
            \State (Device)
                $Z_{i,b} \gets \log(Y_{i,b})$;
                \label{op: log-Y}
            \State (Device)
                $
                    J_{i,b} \gets 
                    e^{\hat{Z} \hat{A}}
                    [ \, B_1^\top \; \cdots \; B_n^\top \,]
                $; \label{op: GEMM}
            \State (Device) Multiply $\diag(e^{-\boldz_i},1,1)$ to the $i$-th block of $J$ for $i=1,\ldots,p$;
        \EndProcedure
    \end{algorithmic}
\end{algorithm}

To hide latency, the evaluations of the extended Jacobian matrices
are carried out in batches.
The $p$ points $\hat{\boldy}_1,\ldots,\hat{\boldy}_p$
are divided into batches of some given size $b$.
For the $i$-th batch, \Cref{alg: evaluation}, evaluates
the extended Jacobian matrices $\hat{D}_{\hat{\boldy}} \homog{\hat{H}} (\hat{\boldy}_k)$
for $k = ib+1, \ldots, ib+b$.
In it, $Z_{i,b}$ and $Y_{i,b}$ are the submatrices of $Z$ and $Y$,
respectively, corresponding to the $i$-th batch of points.
The results are stored in the corresponding submatrix $J_{i,b}$ in $J$.
The evaluation for multiple batches are carried out concurrently,
so that operation~\ref{op: log-Y} of \Cref{alg: evaluation} is applied
to one batch of points while the result of operation~\ref{op: GEMM} 
for the previous batch is being written to device global memory,
thereby effectively hiding latency.
The optimal batch size $b$ is dependent on specific hardware characteristics
(e.g. the relative speed at which the target GPU can carry out 
computational-intensive tasks vs memory-intensive tasks).

\begin{algorithm}
    \caption{Computation of Euler and Newton directions}
    \label{alg: euler-newton}
    \begin{algorithmic}[1]
        \Procedure{EulerNewton}{$i,b$}
            \State (Device) Compute the pointer of $\hat{D}_{\hat{\boldy}} \homog{H}(\boldy_i)$ within $J$;
            \State (Device) 
                Compute the QR decomposition of
                $
                    (\hat{D}_{\hat{\boldy}} \homog{H}(\boldy_i))^\top = 
                    Q_i
                    \begin{bmatrix}
                        R_i \\
                        \boldzero_{2 \times n}
                    \end{bmatrix}
                $;
                \label{op: QR}
            \State (Device) $V_i \gets \begin{bmatrix} \bolde_{N+1} \\ \bolde_{N+2} \end{bmatrix} Q_i^*$;
            \State (Device) Perform row reduction so that rightmost $2 \times 2$ block of $V_i$ is $I$.
        \EndProcedure
    \end{algorithmic}
\end{algorithm}

The Euler and Newton directions are computed by \Cref{alg: euler-newton}.
The results, as $2 \times N$ matrices are stored in the leftmost columns
of $V_i$'s, which reside in device (global) memory.





\section{Preliminary implementation and experiments}

As a proof of concept,
the main algorithms are implemented through \tech{CUDA} framework
and supports \tech{NVIDIA}'s general-purpose GPUs (GPGPUs).%
\footnote{
    Only recent GPUs of ``Pascal'' micro-architecture or newer
    are directly supported due to the use of newer features.
    Tensor core features are only supported on
    \tech{NVIDIA A100} and \tech{A30}
    that have double precision tensor cores.
}
In particular, it relies on the \tech{cuBLAS} library
for efficient usage of GPU resource.

In the base implementation,
the general matrix multiplication (\tech{GEMM}) operations
is carried out by the function \texttt{cublasZgemm},
which performs complex matrix multiplications in double precision.
The QR decomposition computation,
in \Cref{alg: euler-newton} operation~\ref{op: QR},
is carried out by the function \texttt{cublasZgeqrfBatched},
which automatically schedules the decomposition of a large number of matrices
and effectively hiding latency.

In the variant that supports tensor core operations,
\Cref{alg: evaluation}
(specifically operation~\ref{op: GEMM} of \Cref{alg: evaluation})
runs on tensor cores with much greater efficiency.%
\footnote{
    Even though \tech{NVIDIA}'s GPUs for high performance computing
    starting from \tech{V100}, release in 2017, all support tensor cores,
    only the more recent \tech{A100} and \tech{A30}
    directly support double precision computation on tensor cores.
    As a result, tensor core operations are disabled for GPUs
    older than \tech{A100} and \tech{A30}.
}
In that case, since \Cref{alg: euler-newton} runs on \tech{CUDA} cores,
in principle, these two algorithms can run concurrently
on different batches of points.

\Cref{tab: GPU tests} shows tests of this experimental implementation
on the well known ``\tech{cyclic-14}'' system,
which is a system of 14 polynomial equations in 14 (complex) variables
with a total of 184 distinct monomials.
Mobile (notebook computer) GPU \tech{NVIDIA Quadro T2000},
\tech{NVIDIA Tesla K80}%
\footnote{
    \tech{Tesla K80} contains two GPUs,
    only one of them is used in all the testings.
}
(released in 2014),
and \tech{NVIDIA V100} (released in 2017) are used in these testings.
The results represent the average of 3 different runs.
To keep the test results consistent,
a simplified ``Euler-Newton'' step is used,
in which a single Euler prediction step 
is followed by a single Newton iteration.%
\footnote{
    In actual path trackers, the number of Newton iterations
    depends on the quality of the prediction.
    Here, we fix the number of iterations to be exactly 1
    so that the timing information from different runs can be compared.
}
In all runs, the batch size is set to $\frac{1}{4} p$,
where $p$ is the number of points to be processed simultaneously.
This is likely far from optimal for hiding latency.
Further testings and analysis are needed to determine the optimal batch size.
Similarly, \Cref{tab: GPU tests chandra}
shows the test results on the \tech{Chandra-24} system,
which is a system of 24 polynomial equations in 24 (complex) variables
with a total of 324 distinct monomials.
Both tests suggest the very promising scalability of this GPU-accelerated approach.
In particular, the timing information on \tech{NVIDIA V100} GPU
appear to be constant regardless the number of points are processed.
This indicates that the computational time (latency)
is still dominated by host-device memory transfer
and other memory-bound operations,
and thus the algorithm can scale even further if memory transfer can be optimized or reduced.

\begin{table}[h]
    \begin{tabular}{cccc}
        \toprule
        GPU         & \makecell{Mobile GPU\\{}\tech{Quadro T2000}} & \tech{Tesla K80} & \tech{V100} \\
        \midrule
        10   points & 0.0744s & 0.0850s & 0.0222s \\
        50   points & 0.0888s & 0.0991s & 0.0274s \\
        250  points & 0.1910s & 0.1092s & 0.0296s \\
        500  points & 0.3330s & 0.2059s & 0.0313s \\
        1000 points & 0.7035s & 0.4120s & 0.0355s \\
        \bottomrule\\
    \end{tabular}
    \caption{
        Wall-clock time, in seconds, 
        for performing 100 consecutive Euler-Newton steps
        for the \tech{cyclic-14} system
        on groups of points of different sizes
        using various \tech{NVIDIA} GPUs.
    }
    \label{tab: GPU tests}
\end{table}

\begin{table}[h]
    \begin{tabular}{cccc}
        \toprule
        GPU         & \makecell{Mobile GPU\\{}\tech{Quadro T2000}} & \tech{Tesla K80} & \tech{V100} \\
        \midrule
        250  points & 0.8219s & 0.3266s & 0.0728s \\
        500  points & 1.7084s & 0.6461s & 0.0757s \\
        1000 points & -       & 1.1902s & 0.0852s \\
        \bottomrule\\
    \end{tabular}
    \caption{
        Wall-clock time, in seconds, 
        for performing 100 consecutive Euler-Newton steps
        for the \tech{chandra-24} system
        on groups of points of different sizes
        using various \tech{NVIDIA} GPUs.
        ``-'' indicate the test cannot be completed as a single run
        due to GPU memory limitation.
    }
    \label{tab: GPU tests chandra}
\end{table}

For a comparison with existing CPU-based algorithms,
we choose the implementation in Hom4PS-3,%
\footnote{
    Best efforts are made to keep the comparisons fair.
    However, the implementation in Hom4PS-3 is fundamentally different
    as additional steps are carried out to store information
    for higher order predictor-corrector algorithms.
    It appears the time consumed by these additional steps is small in comparison.
}
which adopts the basic algorithm in HOM4PS-2.0.
Moreover, using \tech{Intel MKL} (math kernel library),
the \tech{BLAS}-based algorithm proposed here
can also be easily ported back to modern CPU that support \tech{AVX-512} instructions.%
\footnote{
    \tech{Intel MKL} also supports older CPUs that do not support \tech{AVX-512} instructions.
    However, in those cases, hardware acceleration may have minimum effect.
}
\Cref{tab: GPU vs CPU chandra} shows a four-way comparison of these implementations:
\begin{enumerate}
    \item The proposed GPU-accelerated algorithm;
    \item The proposed algorithm implemented on modern CPUs (restricted to a single core);
    \item The proposed algorithm implemented on modern CPUs (using 8 cores);
    \item The existing Hom4PS-3 (CPU-only) algorithm.
\end{enumerate}

The result shows the clear advantage of the GPU-accelerated algorithm
as the proof-of-concept implementation is already showing vastly superior
performance and scalability.
Interestingly, even when ported back to CPUs (using \tech{Intel MKL}),
the proposed algorithm still appear to be faster than the existing
implementation in Hom4PS-3.

\begin{table}[h]
    \small
    \begin{tabular}{ccccc}
        \toprule
        Algorithm & 
        \makecell{GPU version of the\\{}proposed algorithm\\{}based on \tech{CUDA}} &
        \makecell{CPU version of the\\{}same algorithm\\{}using \tech{BLAS} on 1 core} &
        \makecell{CPU version of the\\{}same algorithm\\{}using \tech{BLAS} on 8 cores} &
        \makecell{\tech{Hom4PS-2/3} algorithm\\{}(CPU only)\\{}running on 8 cores} \\
        \midrule
        Time & 
         0.0852s &
        10.7149s & 
         1.3901s &
         3.3090s \\
        \bottomrule\\
    \end{tabular}
    \caption{
        Wall-clock time, in seconds,
        for performing 100 consecutive Euler-Newton steps
        for the \tech{chandra-24} system
        on a group of 1000 points
        on  GPU (\tech{NVIDIA V100})
        and CPU (\tech{Intel Xeon E5-2686 v4})
        using different algorithms.
    }
    \label{tab: GPU vs CPU chandra}
\end{table}

\section{Concluding remarks}

We proposed an extremely simple approach for efficiently evaluating
a homotopy function together with all its partial derivatives,
in the context the polyhedral homotopy method of Huber and Sturmfels,
that maps particularly well to modern GPUs (and similar hardware architectures).
This is done through a reformulation of the problem of evaluating
Laurent polynomial systems together with their derivatives
into general matrix multiplication (\tech{GEMM}) operations.
Based on this, we also proposed a companion algorithm that unifies
the computation of Euler and Newton directions
which are needed in path tracking.
Together, they allow the entire path tracking process to be carried out by GPUs.

Our experimental implementation shows promising efficiency and scalability
when compared to existing CPU-based algorithms.
Further optimizations are needed to realize the full potential of this approach.
In particular, minimizing host-device communication, optimizing memory layout,
and properly hiding latency appear to be the most important tasks.

In the experimental implementation, no effort was made to analyze and improve
the accuracy, stability, and numerical condition of evaluation
of the homotopy function together with its partial derivatives.
The necessity and design of pre-conditioners for ensuring the
numerical quality of this step will also be a crucial task in future studies.



\bibliographystyle{plain}
\bibliography{library.bib,gpu.bib}

\end{document}